\documentclass{mcom-l2}

\theoremstyle{plain}
\newtheorem{thm}{Theorem}[section]

\newtheorem{heur}[thm]{Heuristic}
\newtheorem*{heur*}{Heuristic}

\theoremstyle{definition}
\newtheorem{defin}[thm]{Definition}
\theoremstyle{remark}

\newtheorem{rem}[thm]{Remark}

\usepackage{algorithm,algorithmic}

\algsetup{indent=2em}

\usepackage{tikz,multirow,url}
\usetikzlibrary{decorations.pathreplacing}

\newcommand{\KK}{{\mathbf K}}

\newcommand{\OO}{{\boldsymbol{\mathcal O}}}
\newcommand{\QQ}{{\mathbf Q}}

\newcommand{\ZZ}{{\mathbf Z}}
\newcommand{\Bcal}{{\mathcal B}}

\newcommand{\Dcal}{{\mathcal D}}

\newcommand{\Lcal}{{\mathcal L}}
\newcommand{\Ncal}{{\mathcal N}}
\newcommand{\Pcal}{{\mathcal P}}

\newcommand{\afrak}{{\boldsymbol{\mathfrak a}}}
\newcommand{\bfrak}{{\boldsymbol{\mathfrak b}}}

\newcommand{\pfrak}{{\boldsymbol{\mathfrak p}}}
\newcommand{\qfrak}{{\boldsymbol{\mathfrak q}}}

\newcommand{\DD}{{|\Delta_\KK|}}

\newcommand{\ie}{{\emph{i.e.,}~}}
\newcommand{\ve}{{\varepsilon}}

\newcommand{\Poly}{{\operatorname{Poly}}}
\newcommand{\Res}{{\operatorname{Res}}}

\newcommand{\gen}[1]{\left\langle {#1} \right\rangle}
\newcommand{\id}[2]{\afrak_{#2}^{(#1)}}

\newcommand{\LL}{\textit{\L}}

\title[Class group computations using small polynomials]{Reducing the complexity for class group computations using small defining polynomials}

\author{Alexandre G\'elin}
\address{Laboratoire de Math\'ematiques de Versailles, UVSQ, CNRS, Universit\'e Paris-Saclay, Versailles, France}
\email{alexandre.gelin@uvsq.fr}

\begin{document}

\begin{abstract}
In this paper, we describe an algorithm that efficiently collect relations in class groups of number fields defined by a small defining polynomial. This conditional improvement consists in testing directly the smoothness of principal ideals generated by small algebraic integers. This strategy leads to an algorithm for computing the class group whose complexity is possibly as low as $L_\DD\left(\frac{1}{3}\right)$. 
\end{abstract}

\maketitle


\section{Introduction}
The ideal class group of a number field is a finite abelian group and its computation is a major task in algorithmic algebraic number theory. The case of quadratic number fields was firstly addressed by Shanks~\cite{Sha69,Sha72}. Thanks to the baby-step--giant-step strategy and under the Generalized Riemann Hypothesis (GRH), he reached an exponential runtime $O(\DD^{\frac 1 5})$, where $\Delta_\KK$ denotes the absolute discriminant of the considered number field. 

Hafner and McCurley~\cite{HMC89} then proposed an algorithm in heuristic subexponential time $L_\DD(\frac 1 2, \sqrt{2})$, but only in the restrictive case of imaginary quadratic number fields.
This $L$-notation is classical when presenting index calculus algorithms with subexponential complexity. Given two constants $\alpha$ and $c$ with $\alpha\in [0,1]$ and $c \geq 0$, $L_N(\alpha,c)$ is used as a shorthand for
\[\exp\left((c + o(1))(\log N)^\alpha (\log \log N)^{1-\alpha} \right),\]
where $o(1)$ tends to $0$ as $N$ tends to infinity. We also encounter the notation $L_N(\alpha)$ when specifying $c$ is undesired.

An extension of this latter algorithm to all number fields was the topic of Buchmann's work~\cite{Buc90}, assuming that the extension degree, arbitrary, is fixed. Then he obtained a heuristic runtime $L_\DD(\frac 1 2, 1.7)$. Finally, Biasse and Fieker improved this algorithm and achieved a subexponential complexity for all number fields, without any restriction on the degree: a complexity $L_\DD(\frac{2}{3}+\varepsilon)$ in the general case\footnote{For an arbitrary small $\varepsilon>0$.} and $L_\DD(\frac{1}{2})$ when the extension degree $n$ satisfies $n\leq (\log\DD)^{3/4-\varepsilon}$. Recently, these complexities were reduced to $L_\DD\left(\frac{2\alpha+1}{5}, o(1)\right)$ for number fields in classes $\Dcal_{n_0,d_0,\alpha,\gamma}$ with $\alpha > 3/4$, and $L_\DD\left(\frac{1}{2},\frac{\omega+1}{2\sqrt{\omega}}\right)$ in the other cases --- the classes $\Dcal$ are defined in~\cite{GJ16} and the complexities come from~\cite{Gel18a}.

In addition, there exists some conditional improvements when the defining polynomial of the number field has good properties --- namely small coefficients. Biasse and Fieker~\cite{BF14} achieved an $L_\DD(a)$ complexity with $a$ possibly as low as $\frac{1}{3}$, and this improvement has been widened in~\cite{GJ16} to a larger set of number fields.

\subsection*{Contribution.} In this paper, we focus on a conditional improvement based on the smallness of the defining polynomial. Though ideal-reduction schemes enforce an $L_\DD\left(\frac{1}{2}\right)$ complexity, the solution of the discrete logarithm problem in finite fields in $L_Q\left(\frac{1}{3}\right)$ suggests that we can reach this value for class group computations too. This is the aim of the \emph{sieving strategy}. We first describe the algorithm and extend the results obtained by Biasse in~\cite{Bia14a}. Then we study its complexity, compare it with the results of~\cite{Gel18a} and exhibit the number fields for which this new strategy offers a better complexity than ideal reductions. In addition, we provide an algorithm for solving the Principal Ideal Problem by using techniques close to the ones used for class group computations.

\subsection*{Outline.} The article is organized as follows. In Section~\ref{sec:motiv} we briefly explain how this sieving strategy may speed up class group computation. Then Section~\ref{sec:algo} is devoted to the description of the relation collection algorithm, while Section~\ref{sec:comp} gives the parameter choices together with the complexity analysis according to the classes~$\Dcal$. Section~\ref{sec:sumup} summarizes where each algorithm --- this one and the one based on ideal reduction --- is better than the other in order to give a new state of the art of class group computation. Finally, the solution of the Principal Ideal Problem based on this method is provided in Section~\ref{sec:PIP}.

\section{Motivation} \label{sec:motiv}

As it is explained in~\cite[Section~2]{Gel18a}, computing class groups and regulators in number fields is essentially based on the index calculus method. Within this strategy, the part that determines the complexity is the relation collection, because the linear-algebra step only leads to an additional constant factor in the exponent --- \ie in the second constant in the $L$-notation. The relation collection step, as its name suggests, consists in searching for many principal ideals that split over the factor base $\Bcal = \left\{\pfrak_1,\dotsc,\pfrak_N\right\}$ composed of all prime ideals of norm below a bound $B>0$:
\[\gen{x} \OO_\KK = \prod \pfrak_i ^{e_i} \qquad \text{ for } x \in \OO_\KK.\]

In the general case, without making any assumption on the number fields, the ideal-reduction strategy performs best and leads to a complexity that is at least $L_\DD\left(\frac{1}{2}\right)$. However, there exist conditional improvements when the number field is defined by a \emph{good} polynomial, that is a polynomial having small height. Indeed, in that case, the $\qfrak$-descent strategy described by Biasse and Fieker in~\cite{BF14} and generalized in~\cite{GJ16} allows a complexity between $L_\DD\left(\frac{1}{3}\right)$ and $L_\DD\left(\frac{1}{2}\right)$ for all number fields of small extension degree.

Our new idea that underlies this article is to generate the relations by testing a lot of \emph{small} principal ideals that are generated by algebraic integers of bounded degree and coefficients. The norms of such elements depend on the two bounds used for the degree and on the coefficients and the height of the defining polynomial. This idea was already used in the Number Field Sieve~\cite{LLMP90}. Enge, Gaudry, and Thom\'e~\cite{EG07,EGT11} extend this method to low-degree curves for solving the discrete logarithm problem over such curves in~$L_{q^g}\left(\frac{1}{3}\right)$, where $q$ is the cardinality of the base field and $g$ the genus of the curve. 

Then, Biasse in~\cite{Bia14a} applies the method in the context of class group computations. His result only addresses very specific number fields $\KK$ defined by a polynomial $T$ such that
\begin{equation} \label{eq:BiaSieve} \big[\KK:\QQ\big] \leq O(\log \DD)^\alpha \quad \mbox{and} \quad \log H(T) \leq O(\log \DD)^{1-\alpha} \end{equation}
for an $\alpha$ in the open interval $\left(\frac{1}{3},\frac{2}{3}\right)$. In so doing, he was able to compute the class group in time $L_\DD\left(\frac{1}{3}\right)$ assuming the Extended Riemann Hypothesis (ERH) and under heuristics. We generalize here the sieving strategy to all number fields, obtaining a complexity possibly as low as $L_\DD\left(\frac{1}{3}\right)$.

This method has also been used by Buchmann, Jacobson, Neis, Theobald, and Weber in~\cite{BJN+99} for practical enhancements. Indeed, the sieving strategy definitely outperforms ideal reduction in practice, especially for small-degree number fields.

The $\qfrak$-descent strategy explained in~\cite{BF14}, where elements with small coefficients are searched in lattices of smaller dimension, is, in a certain sense, another way to use these small algebraic integers. However, our method appears easier to understand and its complexity analysis is streamlined: we are able to provide explicitly the second constant in the $L$-notation, which does not sound that simple for the $\qfrak$-descent. In addition, from a practical point of view, as the $\qfrak$-descent only works in small degree $\left(\alpha \leq \frac{1}{2}\right)$, our algorithm should outperform the $\qfrak$-descent, since it does not require iterations nor lattice-reductions.

\section{Deriving relations by sieving}  \label{sec:algo}

In the following, we make use of the classification presented in~\cite{Gel18a} based on the classes $\Dcal$ introduced in~\cite{GJ16}:

\begin{defin}[{\cite[Definition~3.1]{Gel18a}}]
Let $n_0>1$ be a real parameter arbitrarily close to $1$, $d_0 > 0$, \mbox{$\alpha\in [0,1]$} and $\gamma \geq 1 - \alpha$. The class $\Dcal_{n_0,d_0,\alpha,\gamma}$ is defined as the set of all number fields $\KK$ of discriminant $\Delta_\KK$ that admit a monic defining polynomial $T \in \ZZ[X]$ of degree $n$ that satisfies:
\begin{align} \label{eq:ClassD}
\frac{1}{n_0} \left(\frac{\log \DD }{\log\log \DD}\right)^\alpha \quad \leq \quad &n \quad \leq \quad n_0 \left(\frac{\log \DD }{\log\log \DD}\right)^\alpha \qquad \text{and} \nonumber\\
d = \log H(T) \quad &\leq \quad d_0 (\log \DD)^\gamma (\log\log \DD)^{1-\gamma}.
\end{align}
\end{defin}

For a fixed number field $\KK$ in a class $\Dcal_{n_0,d_0,\alpha,\gamma}$, the value $\alpha \in [0,1]$ corresponds to the extension degree so that it is precisely defined. For the second main parameter $\gamma \geq 1-\alpha$, special care should be taken: sometimes it costs too much to reduce the defining polynomial. This issue is addressed in Section~\ref{sec:sumup}: given a number field defined by a polynomial, we study the optimal strategy for computing the class group depending on the parameters. Is the polynomial reduction necessary? Is it better to use ideal-reduction or sieving?

\begin{rem}
We use the terminology \emph{``sieving strategy''} because it closely corresponds to the way to --- efficiently --- implement it. Theoretically, our algorithm only consists in testing for smoothness a huge arithmetic progression of algebraic integers until we have found sufficiently many relations.
\end{rem}

The description of the algorithm we are going to introduce is clear and the algorithm is easily understandable. Difficulties arise when we need to fix the parameters such as the smoothness bound for the factor base and the bounds that describe the sieving space in order to minimize the complexity. To fix the notation, we consider a number field $\KK = \QQ(\theta)$ of degree~$n$ and let $T$ denote the defining polynomial of which $\theta$ is a root.

Let $B >0$ be the smoothness bound that must be determined. We fix the factor base $\Bcal = \left\{\pfrak_1,\dotsc,\pfrak_N\right\}$ as the set of all prime ideals of $\OO_\KK$ whose norm is below~$B$. From the Landau Prime Ideal Theorem~\cite{Lan03}, we know that its cardinality satisfies \[N = \left|\Bcal\right| = B\big(1+o(1)\big).\]

We describe the sieving space by fixing a bound $t>0$ on the degree, together with a bound $S>0$ on the coefficients. Hence we use all the polynomials of degree at most $t$ with coefficients between $-S$ and $S$. These are $(2S+1)^{t+1}$ polynomials, but only half of them are of interest, as algebraic integers $x$ and $-x$ generate the same ideal. Note that we may also avoid algebraic integers built from a reducible polynomial in $\theta$. Indeed, if $x=x_1\cdot x_2$, then the exponents of a relation produced by $x$ equal the sums of the exponents of relations produced by $x_1$ and $x_2$.

Given an algebraic integer $x= \sum_{i=0}^t a_i \theta^i$ and denoting by $A$ the polynomial $A(X) = \sum a_i X^i$, the norm of the principal ideal $\gen{x}$ is given by \[\Ncal \big(\gen{x}\big) = \Ncal_{\KK/\QQ} \big(x\big) = \Res\left(A,T\right).\]

The bounds for the resultants displayed in~\cite[Theorem~7]{BL10} allow us to provide a bound on the field norm of an element given in standard representation. Thus, thanks to the two bounds $t$ on the degree and $S$ on the coefficients, we can derive an upper bound for the norm of the principal ideal $\gen{x}$:
\begin{equation} \label{eq:BoundRes}
  \Ncal\big(\gen{x}\big) \leq \sqrt{t+1}^{\,n}\, \sqrt{n+1}^{\,t}\, H(T)^t\, S^n.
\end{equation}

We also recalled the two heuristics used in~\cite{Gel18a}, as we also need them.

\begin{heur}[{\cite[Heuristic~4.4]{Gel18a}}] \label{heur:smooth}
The probability $\Pcal(x,y)$ that an ideal of norm bounded by $x$ is \mbox{$y$-smooth} satisfies \[\Pcal(x,y) \geq e^{-u (\log u)(1+o(1))} \quad \text{for} \quad u=\frac{\log x}{\log y}.\]
\end{heur}

\begin{heur}[{\cite[Heuristic~4.7]{Gel18a}}] \label{heur:NbRel+}
There exists $K$ negligible compared with~$|\Bcal|$ such that collecting $K\cdot|\Bcal|$ relations suffices to obtain a relation matrix that generates the whole lattice of relations.
\end{heur}

Assuming Heuristic~\ref{heur:smooth}, the previous bound on the norm offers a lower bound on the probability~$\Pcal$ of $B$-smoothness of any principal ideal $\gen{x}$ belonging to the sieving space. Then the~$(2S+1)^{t+1}$ small ideals lead to $(2S+1)^{t+1} \cdot \Pcal$ relations. Assuming Heuristic~\ref{heur:NbRel+}, collecting $N\big(1+o(1)\big)$ relations suffices to derive the class group. Therefore we want the following relation to be satisfied by our choice of parameters:
\begin{equation} \label{eq:NbRel}
  (2S+1)^{t+1} \cdot \Pcal = N\big(1+o(1)\big).
\end{equation}

\begin{rem}
Note that making use of the weaker Heuristic~\ref{heur:NbRel+}, introduced in~\cite{GJ16}, is essential here. Indeed, the factor base may contain ideals of degree $k > t$, that cannot be part of any relations derived from our settings. Because every ideal whose norm is below the Bach bound has a degree smaller than $\log 12 + 2 \log \log \DD$, we know that sieving on degree-$t$ polynomials suffices for our purposes, which was not the case with the heuristic used before, where the relation matrix must have full rank.
\end{rem}

To evaluate the cost of the sieving phase, we need to know the number of ideals we test for smoothness: it is $(2S+1)^{t+1}$. We explain below that the cost of each smoothness test is always negligible.
Then the overall cost of the sieving phase is given by $(2S+1)^{t+1}\left(1+o(1)\right)$.

As the lowest final complexity is obtained when a balance is reached between the cost of the relation collection and the cost of the linear-algebra phase, we also want that
\begin{equation} \label{eq:Balance} (2S+1)^{t+1} = N^{\omega+1}\big(1+o(1)\big), \end{equation} because the linear algebra cost is in $N^{\omega+1}$ (see~\cite[Proposition~4.1]{BF14}).

Before determining the parameters that minimize the complexity, we give an outline of the strategy in Algorithm~\ref{algo:sieve}.

\begin{algorithm}[h]
  \caption{Deriving relations from small algebraic integers}
  \label{algo:sieve}
  \begin{algorithmic}[1]
    \REQUIRE The factor base $\Bcal$, the degree bound $t$ and the coefficient bound $S$.
    \ENSURE The relations stored.
    \FOR {$d$ \textbf{from} $1$ \textbf{to} $t$}
    \FOR {\textbf{all} $(a_0,\dots,a_d) \in [-S,\dotsc,S]^{d+1}$}
    \STATE Fix $x = \sum a_i \theta^i$ and $\afrak = \gen{x}$
    \STATE Test the $B$-smoothness of $\afrak$
    \IF {$\afrak$ is $B$-smooth}
    \STATE Fix $e_i$ such that $\afrak = \prod \pfrak_i^{e_i}$
    \STATE Store the relation $\gen{x} = \prod \pfrak_i^{e_i}$
    \ENDIF
    \ENDFOR
    \ENDFOR
  \end{algorithmic}
\end{algorithm}

We describe in the subsequent sections how to set the parameters for the factor base and the sieving space to achieve the best complexities. We fix $n_0 >1$, $d_0 > 0$, $\alpha \in [0,1]$ and $\gamma \geq 1-\alpha$ and let $\KK$ be a number field that belongs to $\Dcal_{n_0,d_0,\alpha,\gamma}$. We also assume that we know a \emph{good} defining polynomial $T$ that satisfies \[\log H(T) \quad \leq \quad d_0 (\log \DD)^\gamma (\log \log \DD)^{1-\gamma}.\] Let $\theta$ be a primitive element of $\KK$ that is a root of the defining polynomial $T$. As in the discrete logarithm problem in finite fields, we need to distinguish several cases according to the relative sizes of $\alpha$ and $\gamma$. However, the distinctions between the various cases are not as precise as they are for the DLP: we consider small, medium and large degrees and give the corresponding inequalities involving $\alpha$ and $\gamma$.

\section{Complexity analyses} \label{sec:comp}

\subsection{The case of medium degree}

We begin by the medium case, which we define by $\alpha$ and $\gamma$ being of the same magnitude. This includes $\alpha \approx \gamma \approx \frac{1}{2}$, but covers a much wider range as follows from the analysis below. As already discussed at the beginning of~\cite[Section~3]{Gel18a}, the size of the defining polynomial plays a role in the complexity: we only have the inequality $\gamma \geq 1-\alpha$, so that we have no choice but to keep using both $\alpha$ and $\gamma$.

Given that we hope to find an algorithm with runtime $L_\DD\left(\frac{1}{3}\right)$ and given that $\gamma \geq 1-\alpha$ (thus $\alpha+\gamma \geq 1$), we simply conjecture the existence of an algorithm with runtime $L_\DD\left(\frac{\alpha+\gamma}{3}\right)$ and fix the size of the factor base $\Bcal$ as the set of prime ideals of norm at most \[B = \LL_\DD\left(\frac{\alpha + \gamma}{3},c_b\right),\] with $c_b > 0$ to be determined. The notation $\LL$ is identical as the $L$ introduced earlier, except that we have removed the $o(1)$, in order to consider constants: $\LL_N(\alpha,c) = e^{c\left(\log N\right)^\alpha\left(\log \log N\right)^{1-\alpha}}$.

Thanks to Landau's Prime Ideal Theorem~\cite{Lan03}, we know that $N = \left|\Bcal\right| = L_\DD\left(\frac{\alpha+\gamma}{3},c_b\right)$. The sieving space is chosen to consist in all algebraic integers $x=A(\theta)$, built as polynomials in~$\theta$, that satisfy
\begin{equation} \label{eq:med}
\deg A \leq t=c_t \left(\frac{\log \DD }{\log\log \DD}\right)^{\frac{2}{3}(\alpha+\gamma)-\gamma}\!\mbox{and} \; H(A) \leq S=\LL_\DD\!\left(\frac{2}{3}(\alpha+\gamma)-\alpha,c_s\!\right)\!.
\end{equation}
In particular, $\log H(A) \leq c_s (\log \DD)^{\frac{2}{3}(\alpha+\gamma)-\alpha}(\log\log \DD)^{1-\left(\frac{2}{3}(\alpha+\gamma)-\alpha\right)}$.
So these two quantities~are only well defined for $\frac{2}{3}(\alpha+\gamma)-\gamma \geq 0$ and $\frac{2}{3}(\alpha+\gamma)-\alpha \geq 0$, which defines the bounds of the medium-degree case. 

According to Equation~\eqref{eq:BoundRes}, this choice of parameters enables to bound the norm of every principal ideal $\gen{x}$ in the sieving space by
\begin{equation} \label{eq:mednorm}
\Ncal \big(\gen{x}\big) \leq L_\DD \left(\frac{2}{3}(\alpha+\gamma),n_0c_s+d_0c_t\right).
\end{equation}

We deduce from Heuristic~\ref{heur:smooth} that a principal ideal generated by such an $x$ is $B$-smooth with probability
\[\Pcal \geq L_\DD\left(\frac{\alpha+\gamma}{3},\frac{(\alpha+\gamma)(n_0c_s + d_0c_t)}{3c_b} \right)^{-1}.\] 

The size of the sieving space is given by $(2S+1)^{t+1} = \LL_\DD\left(\frac{\alpha+\gamma}{3},c_sc_t\right)$. As usual, this estimation allows us to estimate the number of relations found by combining the two previous results: the number of collected relations is expected to be \[(2S+1)^{t+1}\cdot \Pcal = L_\DD\left(\frac{\alpha+\gamma}{3},c_sc_t-\frac{(\alpha+\gamma)(n_0c_s + d_0c_t)}{3c_b}\right).\]

With $N = L_\DD\left(\frac{\alpha+\gamma}{3},c_b\right)$ and the assumption of Heuristic~\ref{heur:NbRel+} (see Equation~\eqref{eq:NbRel}), we obtain
\[c_sc_t-\frac{(\alpha+\gamma)(n_0c_s + d_0c_t)}{3c_b} = c_b.\]

Another equation between the various constants stems from the balance between the relation collection and the linear algebra, as stated by Equation~\eqref{eq:Balance}. It boils down to \[c_sc_t = (\omega + 1) c_b.\]

From these two equations, we easily express $c_t$ in the other constants and obtain a deg-2 equation in $c_b$, depending on $c_s$:
$3\omega c_sc_b^2 - d_0(\alpha+\gamma)(\omega +1)c_b - n_0(\alpha+\gamma)c_s^2 = 0$. This expression allows us to infer the shape of $c_b$, which is going to give us the final complexity, depending on~$c_s$:
\[c_b = \frac{d_0(\alpha+\gamma)(\omega +1) + \sqrt{d_0^2(\alpha+\gamma)^2(\omega +1)^2 + 12n_0(\alpha+\gamma)\omega c_s^3}}{6 \omega c_s}.\]

It only remains to minimize this quantity as a function of $c_s$. It follows from a straight analysis that the minimum is achieved for $c_s$ satisfying
$c_s^3 = \frac{2d_0^2(\alpha+\gamma)(\omega + 1)^2}{3 n_0 \omega}$, which leads to \[c_b = \left(\frac{4n_0d_0(\alpha+\gamma)^2(\omega +1)}{9\omega^2}\right)^{\frac{1}{3}}.\]

Consequently, the runtime of our algorithm for computing the class group structure and an approximation of the regulator is
\[L_\DD\left(\frac{\alpha+\gamma}{3},\left(\frac{4n_0d_0(\alpha+\gamma)^2(\omega +1)^4}{9\omega^2}\right)^{\frac{1}{3}}\right).\]

\begin{rem}
The first constant may be as low as $\frac{1}{3}$ if $\gamma$ reaches the lower bound $1-\alpha$, \ie \mbox{$\alpha+\gamma=1$}.
\end{rem}

We also mention that in this case, our second constant is better than the one found by Biasse in~\cite{Bia14a}.

This analysis however only holds when the two quantities $\frac{2}{3}(\alpha+\gamma)-\gamma$ and $\frac{2}{3}(\alpha+\gamma)-\alpha$ are non-negative. These conditions offer the limits of our analysis and can be rewritten as \[\frac{1}{3}\left(\alpha + \gamma\right) \leq \alpha \leq \frac{2}{3}\left(\alpha + \gamma\right) \qquad \Longleftrightarrow \qquad \frac{\gamma}{2} \leq \alpha \leq 2\gamma.\] Therefore, it remains to treat the two complementary cases, when either the size of the defining-polynomial height or the extension degree prevails.

\subsection{The small-degree case: when $2\alpha < \gamma$}

The first extreme case we study is when the size of the defining-polynomial height outweighs the extension degree. It corresponds to the left part of the diagrams displayed in~\cite{Gel18a}, where the $\qfrak$-descent strategy works. In these cases, the extension degree satisfies
\[\alpha < \frac{\gamma}{2} \qquad \Longleftrightarrow \qquad \alpha < \frac{1}{3}\left(\alpha + \gamma\right).\]

We are able to reach a final complexity in $L_\DD\left(\frac{\gamma}{2}\right)$ for the relation collection. As~$\alpha$ is relatively small --- below $\frac{\gamma}{2}$ ---  we know that the defining-polynomial reduction algorithm presented in~\cite{GJ16} runs in time $L_\DD\left(\alpha\right)$, which is strictly less than $L_\DD\left(\frac{\gamma}{2}\right)$. Hence this reduction is always negligible compared with the relation collection, so that it can be considered as a precomputation. According to~\cite[Corollary~3.3]{GJ16}, we can also assume $\gamma \leq 1$.

We fix the size of the factor base $\Bcal$ by considering all the prime ideals having norm below \[B = \LL_\DD\left(\frac{\gamma}{2},c_b\right),\] and we have from Landau's theorem that $N = |\Bcal| = L_\DD\left(\frac{\gamma}{2},c_b\right)$. The sieving space is constructed as before, using all polynomials $A$ that satisfy
\begin{equation} \label{eq:small}
  \deg A \leq t=c_t \quad \mbox{and} \quad H(A) \leq S=\LL_\DD\left(\frac{\gamma}{2},c_s\right).
\end{equation}

These adjustments in the definition are motivated by the desire to minimize the norm size. As the height of the defining polynomial is large, we bound the degree of the algebraic integers to guarantee that the norm stays small.

According to Equation~\eqref{eq:BoundRes}, this choice of parameters enables to bound the norm of every principal ideal $\gen{x}$ in the sieving space by
\begin{equation} \label{eq:smallnorm}
\Ncal \big(\gen{x}\big) \leq L_\DD \left(\gamma, d_0c_t\right).
\end{equation}

Assuming Heuristic~\ref{heur:smooth} allows us to have the following inequality satisfied by the probability for a principal ideal generated by such an $x$ to be $B$-smooth:
\[\Pcal \geq L_\DD \left(\frac{\gamma}{2},\frac{d_0 \gamma c_t}{2 c_b} \right)^{-1}.\] 

As the sieving-space cardinality is $(2S+1)^{t+1} = \LL_\DD\left(\frac{\gamma}{2},c_s(c_t+1)\right)$, we obtain the number of collected relations as before and Equation~\eqref{eq:NbRel} results in $c_s(c_t+1)-\frac{d_0\gamma c_t}{2 c_b} = c_b$. Similarly Equation~\eqref{eq:Balance} leads to $c_s(c_t+1) = (\omega + 1)c_b$. From an identical approach as in the previous section, we find the optimal choices for the constants and conclude that the runtime of our algorithm is 
\[L_\DD\left(\frac{\gamma}{2},\left(\frac{d_0\gamma(\omega +1)^2c_t}{2\omega} \right)^\frac{1}{2}\right).\]

\begin{rem}
  The first constant is always between $\frac{1}{3}$ and $\frac{1}{2}$: the upper bound is a consequence of the precomputation made for finding the minimal-height defining polynomial while the lower one comes from $\gamma > \frac{2}{3}(\alpha+\gamma) \geq \frac{2}{3}$. In the second constant, the factor $c_t$ appears so that the complexity depends on the degree of the polynomials we use for sieving. The minimal value is obtained for $c_t=1$, for a runtime in $L_\DD\left(\frac{\gamma}{2},\left(\frac{d_0\gamma(\omega +1)^2}{2\omega} \right)^\frac{1}{2}\right)$.
\end{rem}

\begin{rem}
A possible alternative for the sieving may be to enlarge the sieving space by allowing larger coefficients --- always below $S'=\LL_\DD\left(\gamma-\alpha,o(1)\right)$ --- and to consider only a random subset of size $L_\DD\left(\frac{\gamma}{2},c_s(c_t+1)\right)$ of the sieving space. Using the bound $S'$ does not affect Equation~\eqref{eq:smallnorm} and the complexity is preserved.
\end{rem}

\subsection{The large-degree case: when $\alpha > 2 \gamma$}

In this last case, the extension degree outweighs the size of the defining-polynomial height. It corresponds to the right part of the diagrams displayed in~\cite{Gel18a}. Here we have to work with the input defining polynomial because finding the minimal one costs too much. As the extension degree is large, we opt for sieving polynomials that have small coefficients and large degrees.

We fix the size of the factor base $\Bcal$ by considering all the prime ideals having norm below \[B = \LL_\DD\left(\frac{\alpha}{2},c_b\right),\] and we have from Landau's theorem that $N = |\Bcal| = L_\DD\left(\frac{\alpha}{2},c_b\right)$. The sieving space is constructed using all polynomials $A$ that satisfy
\begin{equation} \label{eq:big}
  \deg A \leq t=c_t\left(\frac{\log \DD }{\log\log \DD}\right)^{\frac{\alpha}{2}} \quad \mbox{and} \quad H(A) \leq S=\LL_\DD\left(0,c_s\right) = (\log \DD)^{c_s}.
\end{equation}

According to Equation~\eqref{eq:BoundRes}, this choice of parameters enables to bound the norm of every principal ideal $\gen{x}$ in the sieving space by
\begin{equation} \label{eq:bignorm}
\Ncal \big(\gen{x}\big) \leq L_\DD \left(\alpha, n_0\left(c_s+\frac{\alpha}{4}\right)\right).
\end{equation}

We deduce from Equation~\eqref{eq:bignorm} and Heuristic~\ref{heur:smooth} that the probability for a principal ideal generated by such an $x$ to be $B$-smooth satisfies
\[\Pcal \geq L_\DD\left(\frac{\alpha}{2},\frac{n_0\alpha(\alpha+4c_s)}{8c_b} \right)^{-1}.\] 

Finally, an identical analysis enables to find the optimal choice for the constants. The final runtime for our class group algorithm based on sieving strategy satisfies
\[L_\DD\left(\frac{\alpha}{2},\left(\frac{n_0\alpha(\alpha+4c_s)(\omega+1)^2}{8\omega}\right)^\frac{1}{2}\right).\]

\begin{rem}
The first constant is always between $\frac{1}{3}$ and $\frac{1}{2}$ since $\alpha > \frac{2}{3}(\alpha+\gamma) \geq \frac{2}{3}$. In the second constant, the constant $c_s$ appears which can be chosen arbitrarily small. The minimal runtime thus becomes $L_\DD\left(\frac{\alpha}{2},\left(\frac{n_0\alpha^2(\omega+1)^2}{8\omega}\right)^\frac{1}{2}\right)$.
\end{rem}

\begin{rem}
Again, it is possible to enlarge the sieving space by allowing the degree to be larger --- always below $t'=c_t \left(\frac{\log \DD }{\log\log \DD}\right)^{\alpha-\gamma-\ve}$ for $\ve > 0$ arbitrarily small --- and to consider only a random subset of the sieving space of size $L_\DD\left(\frac{\alpha}{2},c_sc_t\right)$. Using the bound $t'$ does not affect Equation~\eqref{eq:bignorm} and the complexity is preserved.
\end{rem}

\section{Conclusion on sieving strategy} \label{sec:sumup}

The complexity analyses we have derived in the previous sections assume that we know a small defining polynomial $T$, that is a witness to the fact that $\KK$ belongs to the class $\Dcal$. We recall that the classes $\Dcal$ satisfy \[\Dcal_{n_0,d_F,\alpha,\gamma_F} \quad \subset \quad \Dcal_{n_0,d_0,\alpha,\gamma_0},\] for $n_0, d_0, d_F > 0$, $0 \leq \alpha \leq 1$ and $1-\alpha \leq \gamma_F < \gamma_0$. To identify the best strategy depending on the inputs, we consider a number field $\KK$ defined by a polynomial $T$ such that
\begin{align*}
  \frac{1}{n_0} \left(\frac{\log \DD }{\log\log \DD}\right)^\alpha \leq \deg T \leq n_0 \left(\frac{\log \DD }{\log\log \DD}\right)^\alpha \; \text{and}\\
  \log H(T) \leq d_0 (\log \DD)^{\gamma_0} (\log\log \DD)^{1-\gamma_0}.
\end{align*}

It is easily verified that $\KK$ belongs to $\Dcal_{n_0,d_0,\alpha,\gamma_0}$. In addition we introduce $\gamma_F$ and $d_F$ so that~$\gamma_F$ is the minimal $\gamma$ such that $\KK \in \Dcal_{n_0,d_F,\alpha,\gamma}$.
Thus we consider two different classes to which $\KK$ belongs, namely $\Dcal_{n_0,d_F,\alpha,\gamma_F}$ and $\Dcal_{n_0,d_0,\alpha,\gamma_0}$; note that \[\KK \in \Dcal_{n_0,d_F,\alpha,\gamma_F} \subset \Dcal_{n_0,d_0,\alpha,\gamma_0}.\]

Given the number field $\KK$ defined by the polynomial $T$ as inputs, we study the different options for computing the class group and give the optimal strategy.
Let us first look at the medium-degree case, where $\frac{\gamma_0}{2} \leq \alpha \leq 2\gamma_0$. Necessarily, we have $\alpha \geq \frac{\alpha+\gamma_0}{3} \geq \frac{1}{3}$.
\begin{itemize}
\item When $\alpha \leq \frac{1}{2}$, as $\gamma_0 \leq 2\alpha$, we have $\frac{\alpha+\gamma_0}{3} \leq \frac{1}{2}$ and sieving is the best strategy.
\item When $\frac{1}{2} < \alpha \leq \frac{3}{4}$, the sieving strategy remains optimal as long as $\frac{\alpha+\gamma_0}{3} \leq \frac{1}{2}$. Indeed, beyond this bound, the ideal-reduction strategy becomes less costly and should be preferred. This happens as soon as $\gamma_0 \geq 1$.
\item Similarly, for $\frac{3}{4} < \alpha \leq 1$, the sieving strategy remains optimal as long as $\frac{\alpha+\gamma_0}{3} \leq \frac{2\alpha+1}{5}$. Above this bound, the ideal-reduction strategy becomes the best option. This happens as soon as $\gamma_0 \geq \frac{4}{5}$.
\end{itemize}

The large-degree case is easier to deal with. Provided that $\alpha > 2\gamma_0$, we know that the sieving strategy results in an algorithm with runtime $L_\DD\left(\frac{\alpha}{2}\right)$, between $L_\DD\left(\frac{1}{3}\right)$ and $L_\DD\left(\frac{1}{2}\right)$, as $\alpha > 2\gamma_0$ implies that $\alpha \geq \frac{2(\alpha+\gamma_0)}{3} \geq \frac{2}{3}$. This is always the best option.

\medskip

The small-degree case is when defining-polynomial reduction plays a role. Indeed, we know that its cost is $L_\DD(\alpha)$ while the sieving strategy runs in time $L_\DD\left(\frac{\gamma}{2}\right)$.
Because $\alpha < \frac{\gamma_0}{2}$, we can always perform this reduction as a precomputation. It allows to find the smallest-height defining polynomial and so the minimal $\gamma_F$. This reduction has two outcomes:
\begin{itemize}
\item If $\frac{\gamma_F}{2} < \alpha$, then the sieving strategy has a complexity in $L_\DD\left(\frac{\alpha+\gamma_F}{3}\right)$, which is negligible compared to the cost of the reduction, so that the final runtime is $L_\DD\left(\alpha\right)$. This can only happens when $\alpha > \frac{1}{3}$, since $\alpha + \gamma_F \geq 1$.
\item If $\frac{\gamma_F}{2} > \alpha$, then the sieving strategy has a complexity that outweighs the cost of the reduction, so that the final runtime is $L_\DD\left(\frac{\gamma_F}{2}\right)$. This value is between $L_\DD\left(\frac{1}{3}\right)$ and $L_\DD\left(\frac{1}{2}\right)$, as the reduction algorithm returns a polynomial such that $\gamma_F \leq 1$ --- this is a direct consequence of~\cite[Corollary~3.3]{GJ16} --- and because $\gamma_F > 2\alpha$ implies that $\gamma_F \geq \frac{2(\alpha+\gamma_0)}{3} \geq \frac{2}{3}$. This is the only option when $\alpha < \frac{1}{3}$. 
\end{itemize}

\begin{figure}[h]
  \centering
  \begin{tikzpicture}[scale=0.9]
    \fill[cyan!50] (0,0) rectangle (12,3);
    \fill[cyan!50] (9,0) -- (9,3) -- (12,3.6) -- (12,0) -- cycle;

    \draw[thick, cyan] (6,3) -- (9,3);
    \draw[thick, cyan] (9,3) -- (12,3.6);
    \draw[thick, cyan] (12,0) -- (12,3.6);

    \fill[cyan!40] (0,3) -- (4,2) -- (8,2) -- (12,3) -- (12,3.6) -- (9,3) -- (6,3) -- cycle;
    \fill[cyan!30] (0,3) -- (4.66,2.33) -- (8,2.33) -- (12,3.2) -- (12,3.6) -- (9,3) -- (6,3) -- cycle;
    \fill[cyan!20] (0,3) -- (5.3,2.66) -- (8,2.66) -- (12,3.4) -- (12,3.6) -- (9,3) -- (6,3) -- cycle;
    \draw[dashed, cyan] (0,3) -- (4.66,2.33) -- (8,2.33) -- (12,3.2);
    \draw[dashed, cyan] (0,3) -- (5.3,2.66) -- (8,2.66) -- (12,3.4);
    \draw[dashed, thick, cyan] (4,2) -- (6,3) -- (0,3) -- (4,2) -- (8,2) -- (12,3);
    \draw[thick, cyan] (6,3) -- (9,3) -- (12,3.6) -- (12,3);

    \node at (2,1) [violet] {$L_\DD\left(\max\left(\alpha,\frac{\gamma_F}{2}\right)\right)$};
    \node at (6,1) [violet]{$L_\DD\left(\frac{\alpha+\gamma_0}{3}\right)$};
    \node at (10,1) [violet]{$L_\DD\left(\frac{\alpha}{2}\right)$};

    \draw[->] (0,0) -- (0,5); 
    \node at (-0.5,5) {$a$};
    \foreach \x/\xtext in {0,2/{\frac{1}{3}},3/{\frac{1}{2}},3.6/{\frac{3}{5}}}
    {\draw (0.1cm,\x) -- (-0.1cm,\x) node[left] {$\xtext\strut$};}
    \draw[->] (0,0) -- (13,0);
    \foreach \x/\xtext in {0,3/{\frac{1}{4}},4/{\frac{1}{3}},6/{\frac{1}{2}},8/{\frac{2}{3}},9/{\frac{3}{4}},12/1}
    {\draw (\x,0.1cm) -- (\x,-0.1cm) node[below] {$\xtext\strut$};}
    \node at (13,-0.5) {$\alpha$};

    \draw[loosely dashed, gray!50] (0,3.6) -- (12,3.6);
    \draw[loosely dashed, cyan] (9,0) -- (9,3);
    \draw[loosely dashed, cyan] (6,0) -- (6,3);
    \draw[loosely dashed, cyan] (4,0) -- (4,2);
    \draw[loosely dashed, cyan] (8,0) -- (8,2);
  \end{tikzpicture}
  \caption{Complexity obtained with our sieving strategy.}
  \label{fig:sieve}
\end{figure}
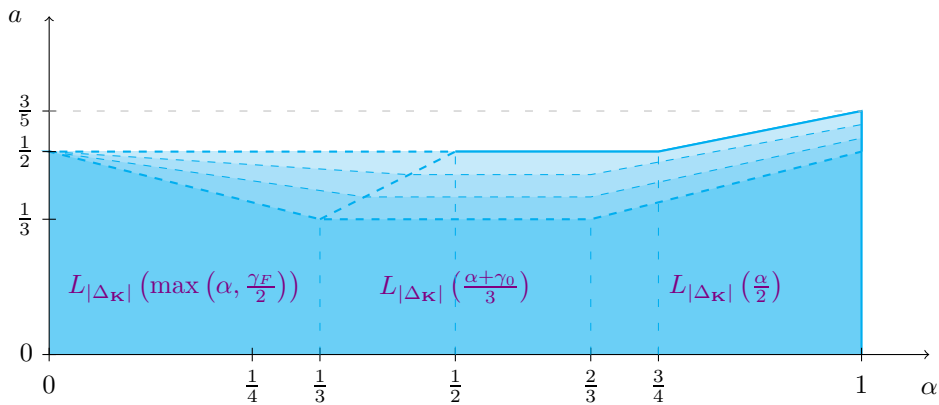

The results of this analysis are summarized in Table~\ref{tab:Choice}. We also give a new diagram for the complexities in Figure~\ref{fig:sieve}.

{\footnotesize 
\begin{table}[ht]
  \centering
  \renewcommand{\arraystretch}{1.5}
  \begin{tabular}{c|ccc}
    {\bf Cond. on $\alpha$} & {\bf Cond. on $\gamma$} & {\bf Strategy} & {\bf Complexity} \\
    \hline
    \multirow{3}{*}{$\alpha \leq \frac 1 2$}
                            & $\gamma_0 \leq 2\alpha$ & Sieving (MD) & $L_\DD \left(\frac{\alpha + \gamma_0}{3} \right)$ \\
    \cline{2-4}
                            & $2\alpha < \gamma_F \leq  \gamma_0$ & Pol. Red. \& Sieving (SD) & $L_\DD \left(\frac{\gamma_F}{2} \right)$ \\
    \cline{2-4}
                            & $\gamma_F < 2\alpha < \gamma_0$ & Pol. Red. \& Sieving (SD) & $L_\DD \left(\alpha\right)$ \\
    \hline
    \multirow{3}{*}{$\alpha > \frac{1}{2}$}
                            & $2\gamma_0 \leq \alpha$ & Sieving (LD) & $L_\DD \left(\frac{\alpha}{2} \right)$ \\
    \cline{2-4}
                            & $\frac{\alpha +\gamma_0}{3} \leq \max\left(\frac{1}{2},\frac{2\alpha+1}{5}\right)$ & Sieving (MD) & $L_\DD \left(\frac{\alpha+\gamma_0}{3}\right)$ \\
    \cline{2-4}
                            & $\frac{\alpha +\gamma_0}{3} > \max\left(\frac{1}{2},\frac{2\alpha+1}{5}\right)$ & Ideal Reduction & $L_\DD \left(\max\left(\frac{1}{2},\frac{2\alpha+1}{5}\right)\right)$ \\
  \end{tabular}
  \caption{Choice of the strategy depending on the input parameters.}
  \label{tab:Choice}
\end{table}}

\section{Application to Principal Ideal Problem} \label{sec:PIP}

In addition to the step forward for class group computations, our results allow us to improve the resolution of another problem: the Principal Ideal Problem (PIP). It consists in finding a generator of an ideal, assuming it is principal. The Short Principal Ideal Problem (SPIP) follows from the PIP by adding the assumption that there exists a \emph{small} generator. The SPIP is the base of several Fully Homomorphic Encryption schemes inspired by the work of Gentry~\cite{Gen09} such as the FHE scheme presented by Smart and Vercauteren at PKC 2010~\cite{SV10} and the multilinear map scheme presented by Garg, Gentry, and Halevi at EuroCrypt in 2013~\cite{GGH13}. Solving the SPIP is a two-stage process that consists of first solving the underlying PIP (on which we focus here), if successful followed by attempts to reduce the generator found to a short one (see~\cite{CDPR16} for instance). Finding a generator of a principal ideal, and even testing the principality of an ideal, are difficult problems in algorithmic number theory, as described in detail in~\cite[Chapter~4]{Coh93} and~\cite[Section~7]{Thi95}.

The general strategy is similar to the one used for the Discrete Logarithm Problem in finite fields. Indeed, for finding the logarithm of an element, two steps are distinguished: first, we find the logarithms of many small elements; second, we express our target element using these small elements and recover its logarithm. It is the same here with our ideal $\afrak$, assumed to be principal. First, we compute the matrix of relations as for class group computations, keeping track of the small elements we have sieved with. Second, we find an ideal~$\bfrak$ that is in the same class as $\afrak$ and that splits over the factor base. Then, linear algebra allows us to recover a generator of $\bfrak$ thanks to the relation matrix and finally, we can solve the~PIP.

\subsection{The descent algorithm}
We first briefly outline the algorithm without fixing the parameters. Indeed, as for class group computations, the optimal parameters choices are derived from the complexity analyses, depending on the number-field exponents $\alpha$ and $\gamma$. In order to bootstrap the descent, we start with a classical BKZ-reduction to obtain an ideal of reasonable norm. Indeed, as the input ideal $\afrak$ is fixed --- the one for which we want a generator --- it can have an arbitrarily large norm. All the ideal reductions are performed on the lattice built from the coefficient embedding~$\varsigma(\afrak)$, as it is described in~\cite[Section~2.2]{BEF+17}. The block-size is fixed so that the complexity of the reduction is strictly below the overall complexity of the algorithm, as it is done in~\cite[Section~5]{Gel18a}. Then the descent consists in a succession of ideal reductions and smoothness tests so that the norms of all ideals involved decrease progressively until they reach the lower bound, given by the smoothness bound used in the class group computations.
We make use of the same result as in~\cite{Gel18a} for the lattice reductions:

\begin{thm}[{\cite[Theorem~4.3]{Gel18a}}] \label{thm:BKZ} 
The smallest vector $v$ output by the BKZ algorithm with block-size~$\beta$ has a norm bounded by
\[\|v\| \quad \leq \quad \beta^{\frac{n-1}{2(\beta-1)}} \cdot (\det \Lcal)^{\frac{1}{n}}.\] The algorithm runs in time $\Poly(n,\log \|B_0\|) \left(\frac{3}{2}\right)^{\beta/2+o(\beta)}$, with $B_0$ the input basis.
\end{thm}

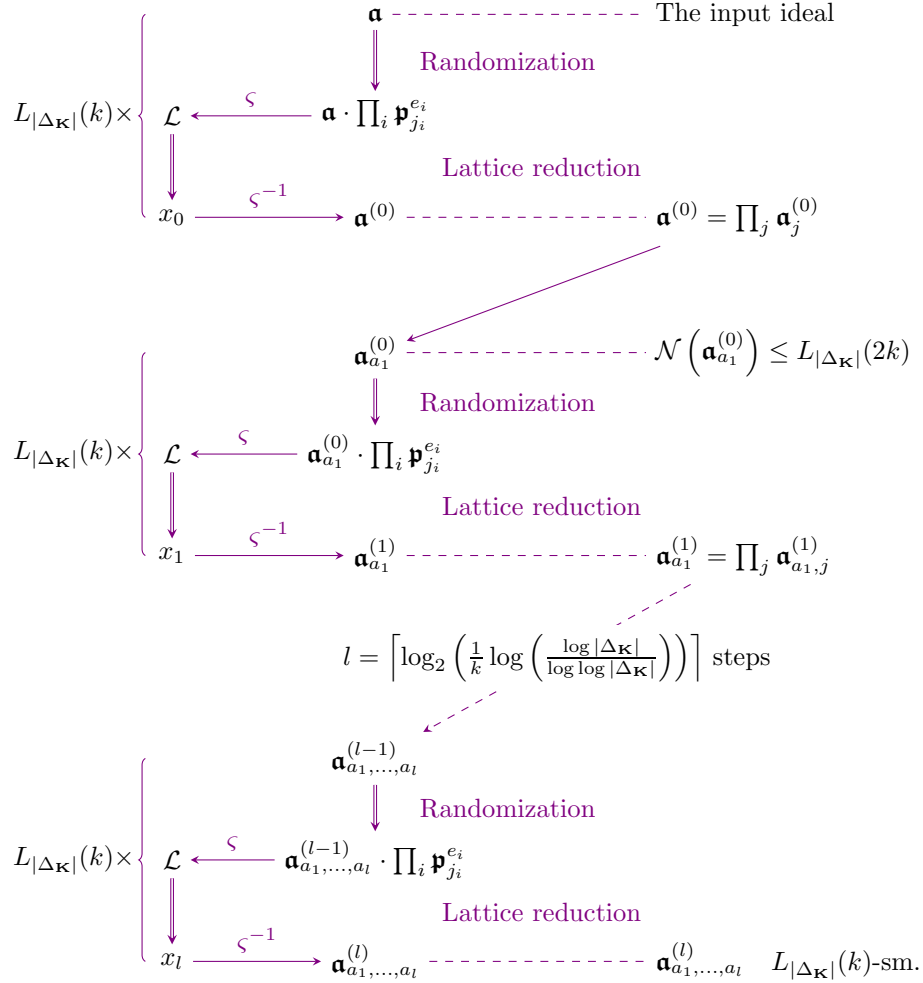
\begin{figure}[h]
  \begin{tikzpicture}[scale=0.9]
    \node (P1) at (3,15) {$\afrak$}; \node[right] (P2) at (7,15) {The input ideal};
    \node (P3) at (0,13.5) {$\Lcal$}; \node (P4) at (3,13.5) {$\afrak \cdot \prod_i \pfrak_{j_i}^{e_i}$};
    \node (P5) at (0,12) {$x_0$}; \node (P6) at (3,12) {$\id{0}{}$}; \node[right] (P7) at (7,12) {$\id{0}{} = \prod_j \id{0}{j}$};
    \draw [violet, decoration={brace,mirror,raise=5pt},decorate] (-0.2,15) -- node[left] {\textcolor{black}{$L_\DD(k) \times\;\;$}} (-0.2,12);
    \draw [violet,>=stealth,->] (P4) -- node[above] {$\varsigma$} (P3); \draw [violet,>=stealth,->] (P5) -- node[above] {$\varsigma^{-1}$} (P6);
    \draw [violet,>=stealth,double,->] (P1) -- node[right] {$\quad$ Randomization} (P4); \draw [violet,>=stealth,double,->] (P3) -- node[right] {\hspace{3cm}$\quad$ Lattice reduction} (P5);
    \draw [violet,dashed] (P1) -- (P2) (P6) -- (P7);

    \node (Q1) at (3,10) {$\id{0}{a_1}$}; \node[right] (Q2) at (7,10) {$\Ncal\left(\id{0}{a_1}\right) \leq L_\DD(2k)$};
    \node (Q3) at (0,8.5) {$\Lcal$}; \node (Q4) at (3,8.5) {$\id{0}{a_1} \cdot \prod_i \pfrak_{j_i}^{e_i}$};
    \node (Q5) at (0,7) {$x_1$}; \node (Q6) at (3,7) {$\id{1}{a_1}$}; \node[right] (Q7) at (7,7) {$\id{1}{a_1} = \prod_j \id{1}{a_1,j}$};
    \draw [violet,decoration={brace,mirror,raise=5pt},decorate] (-0.2,10) -- node[left] {\textcolor{black}{$L_\DD(k) \times\;\;$}} (-0.2,7);
    \draw [violet,>=stealth,->] (Q4) -- node[above] {$\varsigma$} (Q3); \draw [violet,>=stealth,->] (Q5) -- node[above] {$\varsigma^{-1}$} (Q6);
    \draw [violet,>=stealth,double,->] (Q1) -- node[right] {$\quad$ Randomization} (Q4); \draw [violet,>=stealth,double,->] (Q3) -- node[right] {\hspace{3cm}$\quad$ Lattice reduction} (Q5);
    \draw [violet,dashed] (Q1) -- (Q2) (Q6) -- (Q7);

    \node (R1) at (3,4) {$\id{l-1}{a_1,\dotsc,a_l}$};
    \node (R3) at (0,2.5) {$\Lcal$}; \node (R4) at (3,2.5) {$\id{l-1}{a_1,\dotsc,a_l} \cdot \prod_i \pfrak_{j_i}^{e_i}$};
    \node (R5) at (0,1) {$x_l$}; \node (R6) at (3,1) {$\id{l}{a_1,\dotsc,a_l}$}; \node[right] (R7) at (7,1) {$\id{l}{a_1,\dotsc,a_l} \quad L_\DD(k)$-sm.};
    \draw [violet,decoration={brace,mirror,raise=5pt},decorate] (-0.2,4) -- node[left] {\textcolor{black}{$L_\DD(k) \times\;\;$}} (-0.2,1);
    \draw [violet,>=stealth,->] (R4) -- node[above] {$\varsigma$} (R3); \draw [violet,>=stealth,->] (R5) -- node[above] {$\varsigma^{-1}$} (R6);
    \draw [violet,>=stealth,double,->] (R1) -- node[right] {$\quad$ Randomization} (R4); \draw [violet,>=stealth,double,->] (R3) -- node[right] {\hspace{3cm}$\quad$ Lattice reduction} (R5);
    \draw [violet,dashed] (R6) -- (R7);

    \draw [violet,>=stealth,->] (P7) -- (Q1); \draw [violet,dashed,>=stealth,->] (Q7) -- node[fill=white] {\textcolor{black}{$l = \left\lceil \log_2 \left(\frac{1}{k} \log \left(\frac{\log \DD}{\log \log \DD}\right)\right)\right\rceil$ steps}} (R1);
  \end{tikzpicture}
  \caption{The descent algorithm for the medium-degree case.}
  \label{fig:des}
\end{figure}

We now fix the parameters for a degree-$n$ number field $\KK$ that belongs to a class $\Dcal_{n_0,d_0,\alpha,\gamma}$ with $\frac{\gamma}{2} \leq \alpha \leq 2\gamma$. We know that the final complexity is given by $L_\DD\left(\frac{\alpha+\gamma}{3}\right)$, assuming this first constant is small enough --- say below $\frac{1}{2}$. Let us write $k = \frac{\alpha+\gamma}{3}$ for the sake of simplicity. A pattern of the descent is displayed in Figure~\ref{fig:des}.

\subsection*{The initial reduction.}
Let $\afrak$ be the ideal, assumed principal, for which we search for a generator. We may also assume that it is prime, otherwise it suffices to factor it and to work with the prime ideals, which have smaller norms. We can always represent this ideal with its HNF. We obtain an $n \times n$ matrix whose largest coefficient is at most the norm of the ideal~$\Ncal(\afrak)$.

The first reduction consists in performing a BKZ-reduction on the $n$-dimensional lattice~$\varsigma(\afrak)$ with block-size $\beta=(\log \DD)^k$. It permits to exhibit a small vector $v$ that satisfies $\|v\| \leq \beta^{\frac{n-1}{2(\beta-1)}}\Ncal(\afrak)^{\frac{1}{n}}$, as $\det \varsigma(\afrak) = \Ncal(\afrak)$ (see Theorem~\ref{thm:BKZ}). The cost of this lattice reduction is $L_\DD\big(k,o(1)\big)$, provided that the norm $\Ncal(\afrak)$ satisfies $\log \Ncal(\afrak) \leq L_\DD\left(k-\ve\right)$ for $\ve >0$. Therefore, the principal ideal generated by the algebraic integer $x_0 \in \afrak$ corresponding to the vector $v \in \varsigma(\afrak)$ has its norm bounded by $(n+1)^n \cdot H(T)^n \cdot \beta^{\frac{n(n-1)}{2(\beta-1)}} \Ncal(\afrak)$ (using the same technique as for Equation~\eqref{eq:BoundRes}). Finally, denoting by $\id{0}{}$ the unique integral ideal such that $\gen{x_0} = \afrak \cdot \id{0}{}$, we obtain the following upper bound:
\[\Ncal\left(\id{0}{}\right) \leq L_\DD\left(\alpha+\gamma,n_0d_0\right) = L_\DD\left(3k,n_0d_0\right).\]

As we have mentioned, we alternate lattice reductions and smoothness tests. For keeping a complexity in $L_\DD(k)$, we are going to test the ideal $\id{0}{}$ for $L_\DD(2k,s_0)$-smoothness, for $s_0>0$ to be determined. Using ECM algorithm (see~\cite[Appendix~A]{Gel18a}), the cost for a single test is $L_\DD\left(k,\sqrt{2ks_0}\right)$, while the assumption of Heuristic~\ref{heur:smooth} asserts that the probability for $\id{0}{}$ to be $L_\DD(2k,s_0)$-smooth is lower bounded by $L_\DD\left(k,\frac{kn_0d_0}{s_0}\right)^{-1}$. First, this implies that we need to test on average $L_\DD\left(k\right)$ ideals before finding one that is smooth. We then make use of the randomization process used by Biasse and Fieker in~\cite{BF14}. It consists in considering randomized ideals that are products of $\afrak$ with random power-products of small prime ideals --- the ones in the factor base. Clearly, it offers sufficiently many choices for testing $L_\DD\left(k\right)$ ideals. Second, the total runtime for the smoothness tests is given by \[L_\DD\left(k,\frac{kn_0d_0}{s_0}+\sqrt{2ks_0}\right),\] which is minimal for $s_0^3=2k(n_0d_0)^2$, leading to a complexity of \[L_\DD\left(k,\left(\frac{9}{2}k^2n_0d_0\right)^{\frac{1}{3}}\right).\]

\subsection*{Subsequent steps.}
At the beginning of the $i$-th step, we have an ideal $\id{i}{}$ whose norm is upper bounded by $L_\DD\left(k\left(1+\frac{1}{2^i}\right),s_i\right)$. This time, we are going to perform the lattice reduction over a sublattice of $\varsigma\!\left(\id{i}{}\right)$ of dimension $d = c_d \left(\frac{\log \DD}{\log \log \DD}\right)^\delta$, for $0 \leq \delta \leq \alpha$ and $c_d > 0$ to be determined. The reason to look at a sublattice is that it allows to reduce the norms of the ideals that are involved, which is exactly what we want for the descent. 

The BKZ-reduction on this sublattice provides an algebraic integer $x_i \in \id{i}{}$ and so an integral ideal $\id{i+1}{}$ such that $\gen{x_i} = \id{i}{} \cdot \id{i+1}{}$. The upper bound we get on the norm of $\id{i+1}{}$, according to Theorem~\ref{thm:BKZ},~is \[L_\DD(\alpha) \cdot L_\DD\left(\gamma+\delta,d_0c_d\right) \cdot L_\DD(\alpha+\delta-k) \cdot L_\DD\left(\alpha+k\left(1+\frac{1}{2^i}\right)-\delta, \frac{n_0s_i}{c_d}\right).\]

This quantity is minimal when $\gamma+\delta=\alpha+k\left(1+\frac{1}{2^i}\right)-\delta \Longleftrightarrow \delta = \alpha-k\left(1+\frac{1}{2^{i+1}}\right)$ and $c_d^2=\frac{n_0s_i}{d_0}$, which results in the following upper bound for the norm: \[L_\DD\left(k\left(2+\frac{1}{2^{i+1}}\right),2\sqrt{n_0d_0s_i}\right).\]

Again, we want to test this ideal for smoothness and we fix the smoothness bound to $L_\DD\left(k\left(1+\frac{1}{2^{i+1}}\right),s_{i+1}\right)$. This time, the cost for a single ECM is negligible, as given by $L_\DD\left(\frac{k}{2}\left(1+\frac{1}{2^{i+1}}\right)\right)$. The total cost is then inferred from the number of ideals we have to test. Using the same process as for the initial reduction and assuming Heuristic~\ref{heur:smooth}, this number~is \[L_\DD\left(k,\frac{2k\sqrt{n_0d_0s_i}}{s_{i+1}}\right).\]

\subsection*{The final step.}
We fix $l = \left\lceil \log_2 \left(\frac{1}{k} \log \left(\frac{\log \DD}{\log \log \DD}\right)\right) \right\rceil$.
Thus, at step $l$, we have ideals that are $L_\DD\left(k\left(1+\frac{1}{2^l}\right),s_l\right)$-smooth. However, by definition of the $L$-notation,
\begin{multline*}
  \log L_\DD\left(k\left(1+\frac{1}{2^l}\right),s_l\right) \leq s_l (\log \DD)^k (\log \log \DD)^{1-k}\\
  \times \underbrace{\left(\frac{\log \DD}{\log \log \DD}\right)^{1/\log \left(\frac{\log \DD}{\log \log \DD}\right)}}_{=\quad e = \exp(1)} \big(1+o(1)\big),
\end{multline*}
so that we have the inequality $L_\DD\left(k\left(1+\frac{1}{2^l}\right),s_l\right) \leq L_\DD\left(k,e \cdot s_l\right)$.

\begin{rem} \label{rem:e=1}
More precisely, we can go further and get rid of the constant $e$. Indeed, for every $\ve >0$, if $C_\ve$ denotes the smallest integer larger than $\log(1+\ve)^{-1}$, then at step $C_\ve \cdot l$, we only consider ideals that are $L_\DD\left(k,(1+\ve)s_l\right)$-smooth.
\end{rem}

In the end, we want all the ideals involved to have a norm below the smoothness bound we have used for class group computation, \ie 
\begin{equation} \label{eq:SM} e \cdot s_l \leq c_b = \left(\frac{4k^2n_0d_0(\omega+1)}{\omega^2}\right)^{\frac{1}{3}}. \end{equation}
Our approach is to balance the cost of all steps, except the initial one: each one costs $L_\DD\left(k,\left(4k^2n_0d_0y\right)^{\frac{1}{3}}\right)$, for a constant $y>0$ to be determined.
Hence we have, for all $i$, \[\frac{2k\sqrt{n_0d_0s_i}}{s_{i+1}} = 4k^2n_0d_0y \quad \Longleftrightarrow \quad s_{i+1} = \sqrt{s_i} \cdot \left(\frac{4k^2n_0d_0}{y^2}\right)^{\frac{1}{6}}.\]
We deduce that
\begin{align*}
s_l \quad &= \quad s_0^{\frac{1}{2^l}} \cdot \left(\frac{4k^2n_0d_0}{y^2}\right)^{\frac{1}{6} \cdot \left(1+\frac{1}{2}+\cdots+\frac{1}{2^{l-1}}\right)} \\
          &= \quad \left(\frac{s_0 y^{\frac{2}{3}}}{(4k^2n_0d_0)^{\frac{1}{3}}}\right)^{\frac{1}{2^l}} \left(\frac{4k^2n_0d_0}{y^2}\right)^{\frac{1}{3}} \\
          &= \quad \left(\frac{4k^2n_0d_0}{y^2}\right)^{\frac{1}{3}} \big(1+o(1)\big).
\end{align*}

Then, Equation~\eqref{eq:SM} can be rewritten as $\;e \left(\frac{4k^2n_0d_0}{y^2}\right)^{\frac{1}{3}} \leq \left(\frac{4k^2n_0d_0(\omega+1)}{\omega^2}\right)^{\frac{1}{3}}$, \ie $y^2 \geq \frac{e^3\omega^2}{\omega+1}$.
As the number of steps is polynomial in $\log \DD$, the total cost of the $l$ steps of the descent is $L_\DD\left(k,\left(4k^2n_0d_0y\right)^{\frac{1}{3}}\right)$, with $y^2 = \frac{e^3\omega^2}{\omega+1}$.
It outweighs the initial reduction, because $4y > \frac{9}{2}$ for~$\omega \geq 2$.

\begin{rem}
We need to bound the numbers of ideals involved in order to be sure of our final complexity. At each step, we spend time $L_\DD(k)$ for the smoothness tests. It follows that the number of ideals in the decomposition is bounded by $O\left(\left(\frac{\log \DD}{\log\log\DD}\right)^k\right)$. During the descent, the number of ideals is then multiplied by this factor at each step. Finally, the number of ideals at step $l$ is quasi-polynomial $O\left(\left(\frac{\log \DD}{\log\log\DD}\right)^k\right)^l$. In Figure~\ref{fig:des}, indices have been added to the ideals to illustrate this.
\end{rem}

At this point, the only remaining part consists in finding out how to decompose these ideals over the principal ideals collected for building the relation matrix. This is done by solving a linear system $MX=Y$, where $M$ is the relation matrix and $Y$ the valuations vector of the smooth ideal. To be sure that this system has a solution, we need to have a relation matrix of \emph{almost-full} rank. By this unusual term, we only mean that we want all ideals in the factor base involved in the relations, except the ones whose degree is larger than the bound $c_t$. Indeed, they do not appear in a relation because of the parameters we use, but we do not care as they  do not arise either in the descent process --- this is a consequence of the dimensions of the sublattices that we use. The runtime of this part is $L_\DD\left(k,2c_b\right)$ as the matrix of relations is already in HNF.

Finally, we also have $y < \frac{(\omega+1)^4}{\omega^2}$, which means that the complexity for solving the Principal Ideal Problem is the same as the complexity obtained for class group computation. However, we have analyzed the runtime of the descent for the case when the matrix of relations is known.

\begin{rem} \label{rem:comp}
Two improvements can be made to reduce the complexity. First, as explained in Remark~\ref{rem:e=1}, the constant $e$ can be replaced by any other constant larger than and arbitrarily close to 1. Second, if we are only interested in solving the PIP, then the computation of the regulator and the class group structure are useless. Hence, the linear-algebra step boils down to solving a linear system over $\ZZ$, which can be performed in time $L_\DD\left(k,\omega c_b\right)$ using a Las-Vegas algorithm described by Storjohann in~\cite{Sto05}. Then, we can adjust all our parameters replacing $\omega+1$ by~$\omega$. Finally, these enhancements lead to a final complexity for the PIP of \[L_\DD\left(k,\left(\frac{4k^2n_0d_0\omega^4}{(\omega-1)^2}\right)^{\frac{1}{3}}\right).\]
\end{rem}

\begin{rem}
The descent strategy for solving the Principal Ideal Problem is also treated in detail in~\cite{BEF+17}. It is applied in the context of the cryptanalysis of a Fully Homomorphic Encryption Scheme over prime-power cyclotomic fields. The interested reader can find more details there.
\end{rem}

\subsection{The large-degree case} \label{sec:large}

For the present large-degree case, the approach is similar to the previous case, the only difference being the parameters choice. This time, $\alpha > 2\gamma$ and we denote by $k$ the first constant of the class group complexity, \ie $k=\frac{\alpha}{2}$.

We perform the first reduction using a block-size $\beta=c_\beta (\log \DD)^k$. It still costs $L_\DD\big(k,o(1)\big)$ and gives rise to an algebraic integer $x_0$ and an integral ideal $\id{0}{}$ such that $\gen{x_0} = \afrak \cdot \id{0}{}$. The norm of $\id{0}{}$ satisfies \[\Ncal\left(\id{0}{}\right) \leq L_\DD\left(2\alpha-k,\frac{n_0^2}{2c_\beta}\right) = L_\DD\left(3k,\frac{n_0^2}{2c_\beta}\right).\]

We make use of the same randomization process as in the medium-case and obtain a $L_\DD\left(2k,s_0\right)$-smooth ideal in time $L_\DD\left(k,\left(\frac{9k^2n_0^2}{4c_\beta}\right)^{\frac{1}{3}}\right)$, for $s_0^3=\frac{kn_0^4}{2c_\beta^2}$ chosen to minimize this cost. 

The subsequent steps begin with an ideal of norm less than $L_\DD\left(k\left(1+\frac{1}{2^i}\right),s_i\right)$. Then, by fixing $\delta = k\left(1+\frac{1}{2^{i+1}}\right)$, we obtain an ideal $\id{i+1}{}$ such that is norm is upper-bounded by \[L_\DD\left(k\left(2+\frac{1}{2^{i+1}}\right),\frac{n_0s_i}{c_d}\right).\]
so that, assuming Heuristic~\ref{heur:smooth}, we can find a $L_\DD\left(k\left(1+\frac{1}{2^{i+1}}\right),s_{i+1}\right)$-smooth ideal in time \[L_\DD\left(k,\frac{kn_0s_i}{c_ds_{i+1}}\right).\]

In the same way, setting $l = \left\lceil \log_2 \left(\frac{1}{k} \log \left(\frac{\log \DD}{\log \log \DD}\right)\right) \right\rceil$ implies that after step~$l$, the ideals involved are $L_\DD(k,e\cdot s_l)$-smooth; here we want $e \cdot s_l$ to be smaller than~$c_b$.

Let $y>0$ be a constant such that, at each step, the runtime of the smoothness tests is below $L_\DD\left(k,y\right)$. That means that for all $i$, it is the case that $\frac{kn_0s_i}{c_ds_{i+1}} \leq y$. Then, by fixing $c_d = \frac{kn_0}{y} \cdot \left(\frac{es_0}{c_b}\right)^{\frac{1}{l}}$ and $s_{i+1} = s_i \cdot \left(\frac{c_b}{es_0}\right)^{\frac{1}{l}}$, the previous equation is satisfied, resulting in
\[s_l = s_0 \cdot \left(\frac{c_b}{es_0}\right) \qquad \Longleftrightarrow \qquad e \cdot s_l = c_b.\]

As $y>0$ can be chosen arbitrarily small, each step has a runtime in $L_\DD\big(k,o(1)\big)$ and the initial-reduction cost can also be chosen that small, for $c_\beta$ sufficiently large. The remaining part consisting in solving the linear system works in the same way as for the previous case and we can conclude that the complexity of our algorithm for solving the PIP is the same as the complexity of the class group computation. Again, Remark~\ref{rem:comp} holds so that we can reduce the complexity to \[L_\DD\left(k,\left(\frac{k^2n_0\omega^2}{2(\omega-1)}\right)^\frac{1}{2}\right).\]

\subsection{The small-degree case}

Again, we only give a brief summary of the descent. Here we have $2\alpha < \gamma$ and $k$ denotes $\frac{\gamma}{2}$. \bigskip \bigskip

The initial BKZ-reduction provides an ideal of norm below $L_\DD\left(\alpha+\gamma,n_0d_0\right)$ in time $L_\DD\big(k,o(1)\big)$. We can find an ideal that is $L_\DD\left(k+\alpha,s_0\right)$-smooth in time $L_\DD\left(k,\frac{kn_0d_0}{s_0}\right)$ as the cost of a single application of ECM is negligible --- because $\frac{k+\alpha}{2} < k$. 

Then, every subsequent step takes as input an ideal of norm upper bounded by $L_\DD\left(k+\frac{\alpha}{2^i},s_i\right)$. Then, looking for a small vector in the sublattice of dimension $d=c_d \left(\frac{\log \DD}{\log \log \DD}\right)^{\frac{\alpha}{2^{i+1}}}$ leads to a new ideal whose norm is smaller than $L_\DD\left(2k+\frac{\alpha}{2^{i+1}},d_0c_d\right)$. Again we expect, assuming Heuristic~\ref{heur:smooth}, to find one that is $L_\DD\left(k+\frac{\alpha}{2^{i+1}},s_{i+1}\right)$-smooth in time $L_\DD\left(k,\frac{kd_0c_d}{s_{i+1}}\right)$.

At final step $l= \left\lceil \log_2 \left(\frac{1}{\alpha} \log \left(\frac{\log \DD}{\log \log \DD}\right)\right) \right\rceil$, we have $L_\DD(k,e \cdot s_l)$-smooth ideals and we want $e \cdot s_l$ to be smaller than $c_b = \left(\frac{kd_0c_t}{\omega}\right)^{\frac{1}{2}}$. Note that, at this point, $d=c_d e\big(1+o(1)\big)$ for the same reason as above. Hence $c_d$ may be as small as $\frac{1}{e}$ and the cost of the final smoothness test is lower-bounded by \[L_\DD\left(k,\frac{kd_0}{e \cdot s_l}\right) \geq L_\DD\left(k,\frac{kd_0}{c_b}\right) = L_\DD\left(k,\left(\frac{kd_0\omega}{c_t}\right)^{\frac{1}{2}}\right).\]

This last smoothness test dominates the overall complexity of the descent phase, as we can always choose $c_d$ and $x_i$ such that the runtimes of the other smoothness tests become arbitrarily small. In addition, this part is dominated by the class group computation: indeed, $\left(\frac{kd_0\omega}{c_t}\right)^{\frac{1}{2}} \leq (\omega+1)\left(\frac{kd_0c_t}{\omega}\right)^{\frac{1}{2}}$ because $c_t \geq 1 > \frac{w}{w+1}$. Again, we can improve this algorithm as explained in Remark~\ref{rem:comp} and finally get a complexity of \[L_\DD\left(k,\left(\frac{kd_0c_t\omega^2}{\omega-1}\right)^{\frac{1}{2}}\right).\]

\begin{rem}
Thanks to this precise analysis, we are able to derive a precise complexity estimate of the attack presented in~\cite{BEF+17}. Indeed, prime-power cyclotomic fields --- together with their totally real subfields --- asymptotically belong to the class $\Dcal_{1,0,1,1}$. Then the result stated at the very end of Section~\ref{sec:large} implies that the complexity of this attack can be as low as \[L_\DD\left(\frac{1}{2},\frac{\omega}{2\sqrt{2(\omega-1)}}\right).\] Taking $\omega=\log_2 7$, we obtain a runtime for our attack of \[L_\DD\left(\frac{1}{2},0.738\right) = 2^{1.066\cdot n^{\frac{1}{2}}\log n}.\]
\end{rem}

\bibliographystyle{amsalpha}
\bibliography{Biblio}

\end{document}